\documentclass[12pt,a4paper]{article}
\usepackage[english]{babel}
\usepackage{amsthm, amsmath, amssymb}
\usepackage{mathtools}
\usepackage{physics}
\usepackage{xcolor}

\usepackage{float}
\usepackage[none]{hyphenat}
\usepackage{hyperref}
\usepackage{graphicx}
\usepackage{subcaption}
\usepackage{svg}
\usepackage[capitalise,noabbrev,nameinlink]{cleveref}

\makeatletter
\renewcommand{\p@subfigure}{\thefigure.}
\makeatother

\usepackage[a4paper, top=2.5cm, bottom=2.5cm, left=2.5cm, right=2.5cm]{geometry}
\usepackage{csquotes}

\newtheorem{theorem}{Theorem}

\newtheorem{lemma}[theorem]{Lemma}

\newtheorem{proposition}[theorem]{Proposition}
\theoremstyle{definition}
\newtheorem{definition}[theorem]{Definition}
\theoremstyle{definition}

\theoremstyle{definition}

\theoremstyle{definition}

\theoremstyle{definition}

\theoremstyle{definition}

\newcommand{\Q}{\mathbb{Q}}
\newcommand{\R}{\mathbb{R}}

\DeclareMathOperator{\Fix}{Fix}
\DeclareMathOperator{\rfix}{RFix}
\DeclareMathOperator{\sfix}{SFix}
\DeclareMathOperator{\sg}{sign}
\DeclareMathOperator{\Sing}{Sing}
\DeclareMathOperator{\Reg}{Reg}

\DeclareMathOperator{\id}{id}

\title{\textbf{Growth rate of periodic points on closed foliated surfaces} }
\author{\small \textsc{L. Hernández-Corbato, A. Moreno-Becerra, J.J. Sánchez-Gabites}}
\date{}

\begin{document}

\maketitle

\begin{abstract}
  Let $M$ be a closed surface and $f\colon M \to M$ a $\mathcal{C}^{1}$ map with isolated periodic points. If, in addition, $f$ preserves a singular foliation, the number of fixed points of $f^{n}$ is bounded from below, up to an additive constant, by its Lefschetz number. Therefore, $f$ satisfies Shub's conjecture.
\end{abstract}
\begingroup
  \renewcommand{\thefootnote}{}
  \footnotetext{2020 Mathematics Subject Classification. Primary 37C25, 37E99; Secondary 55M20.}
  \footnotetext{Key words and phrases. Growth rate inequality, periodic point, foliation, fixed point index.}
  \footnotetext{The authors are funded by the Spanish Ministerio de Ciencia, Innovación y Universidades through grant PID2024-156663NB-I00.}
\endgroup
In his study of the entropy of diffeomorphisms, M. Shub raised the question of the relationship between the asymptotic growth of the algebraic and geometric numbers of periodic points of $\mathcal{C}^{1}$ self-maps $f\colon M\to M$ of a closed manifold. In particular, Shub stated the following growth rate inequality as Problem 4 (p. 11 of \cite{Shub1974}):
\begin{equation}
  \limsup\frac{1}{n}\log \abs{\Lambda(f^n)} \leq \limsup\frac{1}{n}\log \# \Fix(f^n),
\end{equation}
where $\Lambda$ denotes the Lefschetz number.

If the fixed points of $f^{n}$ are isolated for every $n$, each of them has a well defined fixed point index $i(f^n,p)$ given by the local degree of the map $\id-f^n$ in a sufficiently small neighborhood of $p$. These local indices yield a decomposition of the Lefschetz number via the Lefschetz--Hopf formula:
\begin{equation*}
    \Lambda(f^n) = \sum_{p\in\Fix(f^n)}i(f^n,p).
\end{equation*}
M. Shub and D. Sullivan proved that, for a given $p\in M$, the indices $i(f^{n},p)$ are bounded independently of $n$ if $f$ is $\mathcal{C}^{1}$. Using the preceding formula, they obtained the following corollary  in support of the conjecture \cite{ShubSullivan1974}:
\begin{itemize}
  \item[]  If $f\colon M \to M$ is $\mathcal{C}^{1}$ and the Lefschetz numbers $\Lambda(f^n)$ are not bounded, then the set of periodic points of $f$ is infinite.
\end{itemize}

Motivated by the results of Pugh and Shub \cite{pughshub2014periodic}, Misiurewicz \cite{Misiurewicz2014periodic}, and Graff, Misiurewicz, and Nowak-Przygodzki \cite{graff2019periodic, GraffMisiurewiczNowakPrzygodzki2016}, our work focuses on the case in which $M$ is a compact surface without boundary and $f\colon M \rightarrow M$ is a map of class $\mathcal{C}^{1}$ that preserves a singular foliation. By a singular foliation $\mathcal{F}$ we mean a pair consisting of a discrete subset of $M$, $\Sing(\mathcal{F})$, which we call the set of singular points, and a partition of $M\setminus \Sing(\mathcal{F})$ into disjoint curves that is locally homeomorphic to the collection of horizontal lines in the plane. We call each such curve a leaf and $\Reg(\mathcal{F})\coloneq M\setminus \Sing(\mathcal{F})$ the set of regular points. For any regular point $p$, $L_{p}$ will denote the unique leaf containing it. We can endorse $\Reg(\mathcal{F})$ with the leaf topology, that is, the topology in which a set $A\subset \Reg(\mathcal{F})$ is open if its intersection with each leaf is open in the topology of the leaf given by a parametrization. Every leaf is open by definition, so sets laying on multiple leaves are not connected. Restricted to one leaf, the leaf topology and the topology given by a parametrization coincide, hence the only connected sets are leaf intervals (possibly singletons) and whole leaves.

\begin{definition}
  We say that a map $f\colon M \to M$ preserves a singular foliation if, whenever $p\in M$ and $f(p)$ are regular points, there exists an open interval $I\subset L_p$ containing $p$ whose image is contained in $L_{f(p)}$.
\end{definition}

Observe that this definition is equivalent to saying that, whenever $p$ and $f(p)$ are regular, $f$ is continuous at $p $ in the leaf topology.

\begin{lemma}
  The composition of maps that preserve the foliation also preserves the foliation.
\end{lemma}
\begin{proof}
  Let $f,g\colon M \to M$ be two maps preserving a foliation $\mathcal{F}$, and assume that $p$ and $gf(p)$ are regular. Continuity in the leaf topology establishes the conclusion when $f(p)$ is regular. If $f(p)$ is singular, we can choose a neighborhood $U$ of $f(p)$ containing neither other singularities nor preimages of singularities under $g$. Choosing $I\subset L_{p}$ such that $f(I)\subset U$, for each connected component $J$ of $I\setminus f ^{-1}f(p)$, $f(J)$ is connected in the leaf topology, so contained in a single leaf, and adherent to $f(p)$. Therefore, $gf(J)$ is by the same reasoning connected and adherent to $gf(p)$, thus contained in its leaf.
\end{proof}

We shall prove the following:

\begin{theorem}
  \label{th:maintheorem}
  Let $M$ be a closed surface and $f\colon M \to M $ a $\mathcal{C}^{1}$ map preserving a singular foliation $\mathcal{F}$. Then
    \begin{equation}
      \label{eq:maininequality}
      \abs{\Lambda(f^n)} \leq \# \Fix(f^n) + k,
    \end{equation}
  where $k\geq 0$ is an integer independent of $n$. In particular, Shub's conjecture holds for $f$.
\end{theorem}
The main idea for obtaining this inequality is given by the proposition below, which will be proved later:
\begin{proposition}
  \label{prop:cotaindiceregular}
  Let $f\colon M \to M$ be a map with an isolated fixed point $p$, and suppose that it also preserves a foliation for which $p$ is a regular point. Then $\abs{i(f,p)}\leq 1$.
\end{proposition}

Once this bound for the index at regular points is available, the proof of \cref{th:maintheorem} only requires a way of bounding $i(f^{n},p)$ at singular points. In the $\mathcal{C}^{1}$ case, this bound follows from the Shub-Sullivan theorem.

\begin{proof}[Proof of Theorem \ref{th:maintheorem}]
  For brevity, set $\rfix(f)\coloneq \Fix(f)\cap \Reg(\mathcal{F})$ and $\sfix(f)\coloneq \Fix(f)\cap \Sing(\mathcal{F})$. By the Lefschetz--Hopf formula, the Lefschetz number is the sum of the local indices at each fixed point:
  \begin{equation*}
    \Lambda(f^n) = \sum_{p\in \Fix(f^n)} i(f^n, p) = \sum_{p\in \rfix(f^n)} i(f^n, p)+ \sum_{p\in \sfix(f^{n})} i(f^n, p).
  \end{equation*}
  By a theorem of Shub and Sullivan \cite{ShubSullivan1974}, the sequence of fixed point indices at a point is bounded independently of $ n $: $ \left|i(f^n,p)\right|\leq M(f, p) $. Combining this information with the preceding proposition,
  \begin{align*}
    \left|\Lambda(f^n)\right| \leq \sum_{p\in \rfix(f^{n})} 1 + \sum_{p \in \sfix(f^{n})} M(f,p)\leq  \# \Fix(f^n) + k.
  \end{align*}
\end{proof}

\cref{th:maintheorem} gives a new proof of the exponential growth rate of $\#\Fix(f^{n})$ previously obtained in the literature for some specific foliations \cite{pughshub2014periodic,Misiurewicz2014periodic,GraffMisiurewiczNowakPrzygodzki2016,GraffMisiurewiczNowakPrzygodzki2018, graff2019periodic}, see also \cite{HonoratoEtAl2019GrowthRate}. However, it does not hold in general in the continuous case, as the following examples show.

Consider the foliation of $\mathbb{S}^{2}$ given by the two poles and the parallels. We define $f\colon \mathbb{S}^{2}\to \mathbb{S}^{2}$ to have the north pole as a repeller and the south pole as an attractor, and to map each parallel onto a parallel of strictly lower latitude, traversed twice. Then $\Lambda(f^{n}) = 1+2^{n}$, but the only periodic points of $f$ are the poles.

As another example, define the monopole foliation on $\mathbb{S}^{2}$ by collapsing a meridian arc connecting the two poles in the foliation by parallels just described. The leaf space can be identified with the open interval given by the opposite meridian arc. Let $f\colon \mathbb{S}^{2}\to \mathbb{S}^{2}$ be a map that fixes the unique singularity $p$, acts strictly monotonically on the leaf space, and makes the image of each leaf $L$ traverse its image leaf $L'$ twice, passing through $p$. Then, $\deg(f) = 2$ and $\Lambda(f^{n}) = 1+2^{n}$; however, $p$ is the only periodic point.

If the fixed point indices remain bounded, the conclusion of \cref{th:maintheorem} holds even when the $\mathcal{C}^{1}$ hypothesis fails. An example is given by the following Anosov endomorphism:
\begin{equation*}
    A = \begin{pmatrix}
    3 & 1 \\
    1 & 1
  \end{pmatrix}.
\end{equation*}
This matrix has an eigenvalue $\lambda \notin \Q $ outside the unit disk. The associated linear map induces an endomorphism $T_{A}$ on $\mathbb{T}^{2}$ of degree $d = \det(A) = 2$. The foliation of the plane by straight lines in the unstable direction projects to an irrational foliation of the torus preserved by this endomorphism. Moreover, the action of $T_{A}$ and the foliation are equivariant with respect to the symmetry $z\mapsto -z$, so they induce a map $f$ and a singular foliation $\mathcal{F}$ on the quotient $\mathbb{T}^{2}/\{\pm\id\}$, which is a $2$-sphere. The new foliation $\mathcal{F}$ has four thorn singularities corresponding to the fixed points of the symmetry: three of them are not periodic and the remaining one is fixed, say $p$. The degree of $f$ is $d$, so $\Lambda(f^{n}) = 1+d^{n}$. An analysis of the local dynamics of $f$ at $p$ shows that $i(f^{n},p)=0$ for every $n$, so $k=1$, and we can conclude---although $f$ is not differentiable at the singular points---that the number of fixed points of $f^{n}$ grows at least as $d^{n}$. This conclusion is consistent with the standard computation of $\#\Fix(T_{A}^{n})$.

\begin{proof}[Proof of the Proposition]
    Restricting to a sufficiently small coordinate neighborhood, we may assume that $f$ is defined from an open neighborhood $U$ of the origin to $\R^2$, the origin is its only fixed point, and $f$ preserves the foliation by horizontal lines. Note that $f$ has a coordinate form $f(x,y)=(f_{1}(x,y),f_{2}(y))$.

    We can choose a rectangular neighborhood of the origin $ R\coloneq [-a,a]\times [-b,b]$ such that the homotopy
    \begin{equation*}
    H_{t}(x,y) = ((1-t)f_{1}(x,y)+tf_{1}(x,0),f_{2}(y))
    \end{equation*}
    has no fixed points on the boundary $\partial R$. Indeed, since the origin is the only fixed point of $f$ in $U$, we can choose $a>0$ such that $f_{1}(\pm a,0)\neq \pm a$ and $b>0$ such that $f_{1}(\pm a, y)\neq \pm a$ for every $\abs{y}\leq b$. The only possible fixed points of $H_{t}$ lie in $y = \pm b$. Note that $H_{t}(x,\pm b)=(x,\pm b)$ for some $t$ if and only if $f_{2}(\pm b) =\pm b $ and $x$ lies in between $f_{1}(x,\pm b)$ and $f_{1}(x,0)$. We will show this is not possible when $x\in [0,a]$ and $f_{2}(b) = b$, the remaining cases being analogous. It is sufficient to observe that
    \begin{equation*}
        \sg(f_{1}(x,b)-x) = \sg(f_{1}(a,b)-a) = \sg(f_{1}(a,0)-a) = \sg(f_{1}(x,0)-x),
    \end{equation*}
    where the first and third equalities follow from the fact that the origin is the only fixed point of $f$ and the second follows from the choice of $b$. Consequently, $H_{t}$ has no fixed points on $\partial R$.

    By homotopy invariance and multiplicativity of the index \cite{JezierskiMarzantowicz2006},
    \begin{equation*}
      i(f,0) = i(f, R) = i(H_{0},R) = i(H_{1},R) = i(f_{1}(\cdot,0),[-a,a])\cdot i(f_{2},[-b,b]).
    \end{equation*}
    Since the index of a map from an interval to $\R$ is $ 1$, $-1$ or $0$, we conclude that $\abs{i(f,0)}\leq 1$.
\end{proof}
\medskip
\textbf{AI Declaration.} The authors declare that no artificial intelligence has been used at any stage of the creation of this article.

\bibliographystyle{amsplain}
\bibliography{biblioteca}
\end{document}